\documentclass[preprint]{article}

\usepackage{amssymb}

\title{On the probability of two randomly generated S-permutation matrices to be disjoint}
\author{Krasimir Yordzhev}
\date{}

\newtheorem{theorem}{\bf Theorem}

\newtheorem{corollary}{\bf Corollary}[theorem]

\begin{document}
\maketitle

\begin{center}
Faculty of Mathematics and Natural Sciences\\
South-West University, Blagoevgrad, Bulgaria\\
Email: yordzhev@swu.bg
\end{center}

\begin{abstract}
The concept of S-permutation matrix is considered in this paper. It defines when two binary matrices are disjoint. For an arbitrary $n^2 \times n^2$ S-permutation matrix, a lower band of the number of all disjoint with it S-permutation matrices  is found. A formula for counting a lower band of the number of all  disjoint pairs  of $n^2 \times n^2$  S-permutation matrices is formulated and proven. As a consequence, a lower band of the probability of two randomly generated S-permutation matrices to be disjoint is found.
In particular, a different proof of a known assertion is obtained in the work.
The cases when $n=2$ and $n=3$ are discussed in detail.
\end{abstract}


MSC[2010] 05B20, 60C05

\section{Introduction}\label{intr}

Let $n$ be a positive integer. By $[n]$ we denote the set
$[n] =\left\{ 1,2,\ldots ,n\right\} .$

A \emph{binary} (or   \emph{boolean}, or (0,1)-\emph{matrix})
 is a matrix all of whose elements belong to the set $\mathfrak{B}
=\{ 0,1 \}$. With $\mathfrak{B}_n$ we will denote the set of all  $n
\times n$  binary matrices.

A square binary matrix is called a \emph{permutation matrix}, if there is just one 1 in every row and every column. Let us denote by $\mathcal{P}_n$ the  group of all $n\times n$ permutation matrices, and by $\mathcal{S}_n$ the symmetric group of order $n$, i.e. the group  of all one-to-one mappings of the  set $[n] $ in itself. If $x\in [n]$, $\rho\in \mathcal{S}_n$, then we will denote by $\rho (x)$ the image of the element $x$ in the mapping $\rho$. As is well known,  there is an isomorphism $\theta  : \mathcal{P}_n \rightarrow \mathcal{S}_n$ such that if $A=(a_{ij} )\in \mathcal{P}_n$ and  $\theta (A) =\rho\in \mathcal{S}_n$, then
\begin{equation}\label{theta}
a_{ij} =1 \Leftrightarrow \rho (i) =j.
\end{equation}

Let $i\in [n]$ and let $\rho\in \mathcal{S}_n$. We will call $i$ a \emph{fixed point} or \emph{rencontre} if $\rho (i)=i$. As is well known (see for example  \cite[p. 159]{aigner}) the number $e_p (n)$  of all permutation from $\mathcal{S}_n$  with exactly  $p$ rencontres is equal to
\begin{equation}\label{2prim}
e_p (n)=\frac{n!}{p!} \sum_{k=0}^{n-p} \frac{(-1)^k}{k!}
\end{equation}

It is easy to calculate that
\begin{equation}\label{2sekond}
e_{n-1} (n)=0
\end{equation}

When $p=0$, e.i. if  $\rho (i)\ne i$ for all $i=1,2,\ldots ,n$,  $\rho \in \mathcal{S}_n$  is called \emph{derangement}. The number $d_n$ of all derangements in $\mathcal{S}_n$ is equal to
\begin{equation}\label{defrage}
d_n =n! \sum_{k=0}^n \frac{(-1)^k}{k!}
\end{equation}

Let $n$ be a positive integer and let $A\in \mathfrak{B}_{n^2}$ be a $n^2 \times n^2$ binary matrix.  With the help of $n - 1$ horizontal lines and $n - 1$ vertical lines $A$ has been separated into $n^2$ of number non-intersecting $n\times n$ square sub-matrices $A_{kl}$, $1\le k,l\le n$, e.i.
\begin{equation}\label{matrA}
A =
\left[
\begin{array}{cccc}
A_{11} & A_{12} & \cdots & A_{1n} \\
A_{21} & A_{22} & \cdots & A_{2n} \\
\vdots & \vdots & \ddots & \vdots \\
A_{n1} & A_{n2} & \cdots & A_{nn}
\end{array}
\right] .
\end{equation}

The sub-matrices $A_{kl}$, $1\le k,l\le n$ will be called \emph{blocks}.

 A matrix $A\in \mathfrak{B}_{n^2}$ is called an \emph{S-permutation} if in each row, each column, and each block of $A$ there is exactly one 1. Let the set of all $n^2 \times n^2$ S-permutation matrices be denoted by $\Sigma_{n^2}$.

As it is proved in \cite{dahl} the cardinality of the set of all S-permutation matrices is equal to
\begin{equation}\label{broy}
\left| \Sigma_{n^2} \right| =  \left( n! \right)^{2n} .
\end{equation}

Two binary   matrices $A=(a_{ij} )\in \mathfrak{B}_{m}$ and $B=( b_{ij} )\in \mathfrak{B}_{m}$ will be called \emph{disjoint} if there are not elements with one and the same indices $a_{ij}$ and $b_{ij}$ such that $a_{ij} =b_{ij} =1$, i.e. if $a_{ij} =1$ then $b_{ij} =0$ and if $b_{ij} =1$ then $a_{ij} =0$, $1\le i,j\le m$.

The concept of S-permutation matrix was introduced by Geir Dahl  \cite{dahl} in relation to the popular Sudoku puzzle.

 Obviously a square $n^2 \times n^2$ matrix $M$ with elements of $[n^2 ] =\{ 1,2,\ldots ,n^2 \}$ is a Sudoku matrix if and only if there are  matrices $A_1 ,A_2 ,\ldots ,A_{n^2} \in\Sigma_{n^2}$, each two of them are disjoint and such that $P$ can be given in the following way:
\begin{equation}\label{disj}
M=1\cdot A_1 +2\cdot A_2 +\cdots +n^2 \cdot A_{n^2}
\end{equation}

In \cite{Fontana} Roberto Fontana offers an algorithm which randomly gets a family of
$n^2 \times n^2$ mutually disjoint S-permutation matrices, where $n = 2, 3$. In $n = 3$ he
ran the algorithm 1000 times and found 105 different families of nine mutually
disjoint S-permutation matrices. Then using (\ref{disj}) he obtained $9! \cdot 105 = 38\; 102\; 400$ Sudoku matrices.
In relation with Fontana's algorithm, it looks useful to calculate
the probability of two randomly generated S-permutation matrices to be disjoint.

For the classification  of all non defined concepts and notations
as well as for common assertions which have not been proved here,
we recommend sources \cite{aigner,petrov,Sachkov}.

\section{Main results}

\begin{theorem}\label{th1}
Let $A\in \Sigma_{n^2}$. Then the cardinality $\xi_{n}$ of the set of all matrices  $B\in \Sigma_{n^2}$  which are disjoint with $A$ is equal to
\begin{equation}\label{main}
\xi_{n} =\left( n! \right)^{2n} \left(  \sum_{k=0}^n \frac{(-1)^k}{k!} \right)^n \left( 2-  \sum_{k=0}^n \frac{(-1)^k}{k!} \right)^n +R_n ,
\end{equation}
where $R_n \ge 0 .$
\end{theorem}

Proof: It is easy to see that if $A,B\in\Sigma_{n^2}$ then there are only $n^2 \times n^2$  permutation  matrices $C,D\in \mathcal{P}_{n^2}$ of the type
\begin{equation}\label{CD}
C= \left[
\begin{array}{ccccc}
C_{1} & O & O & \cdots & O \\
O & C_{2} & O & \cdots & O \\
O & O & C_3 & \cdots & O \\
\vdots & \vdots & \vdots & \ddots & \vdots \\
O & O & O &  \cdots & C_{n}
\end{array}
\right] ,\quad
D= \left[
\begin{array}{ccccc}
D_{1} & O & O & \cdots & O \\
O & D_{2} & O & \cdots & O \\
O & O & D_3 & \cdots & O \\
\vdots & \vdots & \vdots & \ddots & \vdots \\
O & O & O &  \cdots & D_{n}
\end{array}
\right]
\end{equation}
where $C_i ,D_i \in \mathcal{P}_n$ are permutation $n\times n$ matrices, $i=1,2,\ldots , n$ and $O$ is the zero $n\times n$ matrix and such that
\begin{equation}\label{B=CAD}
B=CAD .
\end{equation}

Let for all $i=1,2,\ldots ,n$ the elements $\theta (C_i )\in \mathcal{S}_n$ be derangements, or for all $i=1,2,\ldots ,n$ the elements $\theta (D_i )\in \mathcal{S}_n$ be derangements, where $\theta  : \mathcal{P}_n \rightarrow \mathcal{S}_n$  is the isomorphism  defined by formula (\ref{theta}).  It is easily seen that in this case, which we will call \emph{basic case} the S-permutation matrices $A$ and $B=CAD$ are disjoint. This, according to (\ref{defrage}), can be done in
$$\nu_{n} =\left( d_n \cdot n! +n!\cdot d_n -d_n^2 \right)^n =\left[ 2\left( n! \right)^2 \sum_{k=0}^n \frac{(-1)^k}{k!} -(n!)^2 \left( \sum_{k=0}^n \frac{(-1)^k}{k!} \right)^2 \right]^n =$$
$$ =\left( n! \right)^{2n} \left(  \sum_{k=0}^n \frac{(-1)^k}{k!} \right)^n \left( 2-  \sum_{k=0}^n \frac{(-1)^k}{k!} \right)^n $$
ways. Obviously $\xi_n =\nu_n +R_n,$ where $R_n \ge 0$ and $ R_n $ counts all cases when there are $i,j \in [n]$, such that $\theta (C_i )$ and $\theta (D_j )$ are not derangements, but $A$ and $B=CAD$ are not disjoint.

\hfill $\Box$

\begin{corollary}\label{ccrrll1}
Formula (\ref{broy}) immediately follows from (\ref{B=CAD}).

\hfill $\Box$
\end{corollary}

Corollary \ref{ccrrll1} is a different proof of the proposition 1 of \cite{dahl}.

\begin{corollary}\label{th2}
The cardinality $\eta_{n}$ of the set of all  disjoint non-ordered  pairs of $n^2 \times n^2$ S-permutation matrices is equal to
$$\eta_{n} =\frac{(n!)^{2n}}{2} \xi_n =$$
$$=\frac{(n!)^{2n}}{2}\left[ \left( n! \right)^{2n} \left(  \sum_{k=0}^n \frac{(-1)^k}{k!} \right)^n \left( 2-  \sum_{k=0}^n \frac{(-1)^k}{k!} \right)^n +R_n \right],$$
where
$$R_n \ge 0 .$$
\end{corollary}

The proof follows directly from Theorem \ref{th1}, formula (\ref{broy}) and having in mind that the ''disjoint'' relation is symmetric and antireflexive.

\hfill $\Box$

\begin{corollary}\label{th3}
The probability $p(n)$ of two randomly generated $n^2 \times n^2$ S-permutation matrices to be disjoint is equal to
$$p(n) = \frac{\displaystyle \xi_n}{\displaystyle  \left( n! \right)^{2n} -1} =\frac{\displaystyle \left( n! \right)^{2n} \left(  \sum_{k=0}^n \frac{(-1)^k}{k!} \right)^n \left( 2-  \sum_{k=0}^n \frac{(-1)^k}{k!} \right)^n +R_n}{\displaystyle  \left( n! \right)^{2n} -1} .$$
\end{corollary}

Proof: Applying Corollary \ref{th2} and formula (\ref{broy}), we obtain:

$$p(n)= \frac{\displaystyle \eta_{n}}{\displaystyle {\left| \Sigma_{n^2} \right| \choose 2}} =$$
$$= \frac{\displaystyle \frac{(n!)^{2n}}{2}\left[ \left( n! \right)^{2n} \left(  \sum_{k=0}^n \frac{(-1)^k}{k!} \right)^n \left( 2-  \sum_{k=0}^n \frac{(-1)^k}{k!} \right)^n +R_n \right]}{\displaystyle \frac{\left( n! \right)^{2n} \left( \left( n! \right)^{2n} -1\right) }{2}} =$$

$$ = \frac{\displaystyle \left( n! \right)^{2n} \left(  \sum_{k=0}^n \frac{(-1)^k}{k!} \right)^n \left( 2-  \sum_{k=0}^n \frac{(-1)^k}{k!} \right)^n +R_n}{\displaystyle  \left( n! \right)^{2n} -1} =\frac{\displaystyle \xi_n}{\displaystyle  \left( n! \right)^{2n} -1} .$$

\hfill $\Box$

\section{The application of Theorem \ref{th1} and its corollaries to $n=2$ and $n=3$}\label{sec3}\label{application}
We will find the value of  $\xi_n$, $\eta_n$ and $p(n)$   when $n = 2$ and $n=3$. For larger values ​​of $n$ additional efforts are required .

\subsection{Consider $n=2$}

Considering (\ref{2sekond}), when  $n=2$, each of the matrices $C_1$, $C_2$, $D_1$ and $D_2$ defined by (\ref{CD}) and (\ref{B=CAD}) is either derangement or the identity matrix $E_2$, where $\theta  : \mathcal{P}_n \rightarrow \mathcal{S}_n$ is the isomorphism defined by formula (\ref{theta}).

Let $C_1 =E_2$. Then a necessary condition for the matrices  $A$ and $B = CAD$ to be disjoint is $\theta (D_1 )$ and $\theta (D_2 )$ to be derangements. This, however, is the main case.  For a definition of this term see in the proof of Theorem  \ref{th1}. Similarly, we consider the cases $C_2 =E_2$, $D_1 =E_2$ or $D_2 =E_2$. Therefore, if $n=2$ is satisfied $R_2 =0.$ So we get:
$$\xi_{2} = 2^{4} \left(  \sum_{k=0}^2 \frac{(-1)^k}{k!} \right)^2 \left( 2-  \sum_{k=0}^2 \frac{(-1)^k}{k!} \right)^2 =9$$
$$\eta_{2} =\frac{2^{4}}{2} \xi_2 =72$$
$$p(2) = \frac{\displaystyle \xi_2}{\displaystyle   2^{4} -1} =\frac{3}{5}$$

The accuracy of the obtained results is confirmed by the value of $\eta_2 =72$ obtained with other methods in  \cite{ISRN}, where graph theory techniques have been used \cite{Yordzhev2013}.

\subsection{Consider $n=3$}

For the calculation of $R_3$, it is necessary to consider all the possibilities of permutation matrices $C_i$ and $D_j$, $1\le i,j\le 3$, defined with the help of (\ref{CD}) and (\ref{B=CAD}). Let $\mathcal{E}_0$ denote the set of all derangements in $\mathcal{S}_3$, i.e. $\mathcal{E}_0 =\{ \rho\in \mathcal{S}_3 \; | \: \rho (i)\ne i,\; i=1,2,3\}$. It is necessary to consider the following cases:

i) There are only $i,j\in\{ 1,2,3\}$ such that $\theta (C_i ),\theta (D_j ) \notin \mathcal{E}_0$, where  $A$  and $B=CAD$ are disjoint.

ii) There are only $i,j\in\{ 1,2,3\}$ such that $\theta (C_i ),\theta (D_j )\in \mathcal{E}_0$, where $A$  and $B=CAD$ are disjoint.

iii) There are $i_1 ,i_2 ,i_3 ,j_1 ,j_2 ,j_3 \in \{ 1,2,3\}$, $i_s \ne i_t$, $j_s \ne j_t$ when  $s\ne t$, such that $\theta (C_{i_1} ),\theta (C_{i_2} ),\theta (D_{j_1} ) \notin \mathcal{E}_0$, $\theta (C_{i_3} ),\theta (D_{j_2} ),\theta (D_{j_3} ) \in \mathcal{E}_0$, where $A$ and $B=CAD$ are disjoint.

iv) There are $i_1 ,i_2 ,i_3 ,j_1 ,j_2 ,j_3 \in \{ 1,2,3\}$, $i_s \ne i_t$, $j_s \ne j_t$ when $s\ne t$, such that $\theta (C_{i_1} ),\theta (D_{j_1} ),\theta (D_{j_2} ) \notin \mathcal{E}_0$, $\theta (C_{i_2} ),\theta (C_{i_3} ),\theta (D_{j_3} )\in \mathcal{E}_0$, where $A$ and $B=CAD$ are disjoint. This case is considered analogous to the case iii.

v) $\theta (C_1 ),\theta (C_2 ),\theta (C_3 )\notin \mathcal{E}_0$ and there is only  $j\in\{ 1,2,3\}$ such that  $\theta (D_j ) \notin \mathcal{E}_0$, where $A$ and $B=CAD$ are disjoint..

vi) $\theta (D_1 ),\theta (D_2 ),\theta (D_3 )\notin \mathcal{E}_0$ and there is only $i\in\{ 1,2,3\}$ such that $\theta (C_i ) \notin \mathcal{E}_0$, where $A$ and $B=CAD$ are disjoint. This case is considered analogous to the case v.

vii) There is only $i\in \{ 1,2,3 \}$ such that $\theta (C_i )\in \mathcal{E}_0$, where $A$ and $B=CAD$ are disjoint.

viii) There is only $j\in \{ 1,2,3 \}$ such that $\theta (D_j )\in \mathcal{E}_0$, where $A$ and $B=CAD$ are disjoint.  This case is considered analogous to the case vii.

After routine calculations, we obtained that $R_3 =19 \; 008$. For $\nu_3 $ we obtain
$$ \nu_3 =\left( 3! \right)^{6} \left(  \sum_{k=0}^3 \frac{(-1)^k}{k!} \right)^3 \left( 2-  \sum_{k=0}^3 \frac{(-1)^k}{k!} \right)^3 =8\; 000,$$
from where
$$\xi_3 =\nu_3 +R_3 =27\; 008.$$
$$\eta_{3} =\frac{(3!)^{6}}{2} \xi_3 =630\; 042\; 624$$
$$p(3) = \frac{\xi_3}{(3!)^{6} -1} =\frac{27\; 008}{46655}\approx 0,579$$

The accuracy of the  obtained results is confirmed by the value of $\eta_3 =630\; 042\; 624$ obtained with other methods in  \cite{ISRN}, where graph theory techniques have been used \cite{Yordzhev2013}.

\section{Conclusions and future work}
As already mentioned in section \ref{application}, the calculation  $R_n$, $\xi_n$, $\eta_n$ and $p(n)$ for $n\ge 4$ still remains to be done.

We do not know a general formula for finding the number of all $n^2 \times n^2$ Sudoku
matrices for each natural number $n\ge 2$ and we consider that this is an open combinatorial
problem. Only some special cases are known. For example, in $n=2$ it is known \cite{Fontana,ISRN} that $$\sigma_2 = 288.$$
In  \cite{Felgenhauer} it has been shown that  there are exactly
$$\sigma_3 = 9! \cdot 72^2 \cdot 2^7 \cdot 27\; 704\; 267\; 971 = 6\; 670\; 903\; 752\; 021\; 072\; 936\; 960 $$
number of $9\times 9$ Sudoku matrices.

The combinatorial problem  to find the number $\mu (n,k)$ of all $k$-tuples $(2\le k\le n)$ mutually disjoint $n^2 \times n^2$ S-permutation matrices is still open for science. The enumerating all disjoint non-ordered  pairs of $n^2 \times n^2$ S-permutation matrices (see Corollary \ref{th2})  brings us closer to its solution.

From (\ref{disj}) it follows that the  equality $\sigma_n = n! \,\mu (n,n)$ is valid. Then we receive:
$$\mu (2,2) =\frac{\sigma_2}{4!} =\frac{288}{24} =12$$
$$\mu (3,3) =\frac{\sigma_3}{9!} =\frac{6\; 670\; 903\; 752\; 021\; 072\; 936\; 960}{362\; 880} =18\; 383\; 222\; 420\; 692\; 992$$

\section{Acknowledgments}

The author wishes to thank the anonymous referees for their valuable suggestions.


\begin{thebibliography}{1}
\expandafter\ifx\csname url\endcsname\relax
  \def\url#1{\texttt{#1}}\fi
\expandafter\ifx\csname urlprefix\endcsname\relax\def\urlprefix{URL }\fi
\expandafter\ifx\csname href\endcsname\relax
  \def\href#1#2{#2} \def\path#1{#1}\fi

\bibitem{aigner}
M.~Aigner, Combinatorial Theory, Classics in Mathematics, Springer-Verlag,
  1979.

\bibitem{dahl}
G.~Dahl, Permutation matrices related to sudoku, Linear Algebra and its
  Applications 430~(8–9) (2009) 2457 -- 2463.
\newblock \href {http://dx.doi.org/10.1016/j.laa.2008.12.023}
  {\path{doi:10.1016/j.laa.2008.12.023}}.

\bibitem{Fontana}
R.~Fontana, Fractions of permutations. an application to sudoku, Journal of
  Statistical Planning and Inference 141~(12) (2011) 3697 -- 3704.
\newblock \href {http://dx.doi.org/10.1016/j.jspi.2011.06.001}
  {\path{doi:10.1016/j.jspi.2011.06.001}}.

\bibitem{petrov}
N.~Petrov, Probability - Teory as a Systemic View of Nature and Sosiety, Trakia
  Universyty, St. Zagora, 2009.

\bibitem{Sachkov}
V.~N. Sachkov, V.~E. Tarakanov, Combinatorics of Nonnegative Matrices,
  Translations of Mathematical Monographs, American Mathematical Society, 2002.

\bibitem{ISRN}
K.~Yordzhev, Bipartite graphs related to mutually disjoint s-permutation
  matrices, ISRN Discrete Mathematics 2012~(Article ID 384068) (2012) 18.
\newblock \href {http://dx.doi.org/10.5402/2012/384068}
  {\path{doi:10.5402/2012/384068}}.

\bibitem{Yordzhev2013}
K.~Yordzhev, On the number of disjoint pairs of s-permutation matrices,
  Discrete Applied Mathematics~(In Press, Corrected Proof).
\newblock \href {http://dx.doi.org/10.1016/j.dam.2013.06.007}
  {\path{doi:10.1016/j.dam.2013.06.007}}.

\bibitem{Felgenhauer}
B.~Felgenhauer, F.~Jarvis,
  \href{http://www.afjarvis.staff.shef.ac.uk/sudoku/sudoku.pdf}{Enumerating
  possible sudoku grids} (2005).
\newline\urlprefix\url{http://www.afjarvis.staff.shef.ac.uk/sudoku/sudoku.pdf}

\end{thebibliography}
\end{document}